\documentclass[12pt]{amsart}
\def\cal{\mathcal}

\usepackage{latexsym}
\usepackage[dvips]{graphics}
\usepackage[dvips]{color}

\font\Bbb=msbm10 scaled 1200

\def\BB#1{\mbox{\Bbb #1}}              
\def\Mm#1{\mbox{\boldmath$\scriptstyle #1$\unboldmath}}
\def\MM#1{\mbox{\boldmath$#1$\unboldmath}}

\def\I{\CC{i}}

\def\C#1{\multicolumn{1}{c}{#1}}                
\def\captionsize{\@setsize\captionsize{10}\ixpt\@ixpt}%

\def\ps@myheadings{\let\@mkboth\@gobbletwo
 \def\@oddhead{\hfill \it \rightmark \hfill \rm \thepage}
 \def\@oddfoot{}
 \def\@evenhead{\rm \thepage \hfill \it \leftmark \hfill}
 \def\@evenfoot{}
 \def\sectionmark##1{}
 \def\subsectionmark##1{}
}

\newenvironment{Eqnarray*}
{\arraycolsep=1.4pt
  \begin{eqnarray*}}
  {\end{eqnarray*}
    \hspace*{-4pt}}

\catcode`\@=11

\font\tenmsx=msam10
\font\sevenmsx=msam7
\font\fivemsx=msam5
\font\tenmsy=msbm10
\font\sevenmsy=msbm7
\font\fivemsy=msbm5
\newfam\msxfam
\newfam\msyfam
\textfont\msxfam=\tenmsx  \scriptfont\msxfam=\sevenmsx
  \scriptscriptfont\msxfam=\fivemsx
\textfont\msyfam=\tenmsy  \scriptfont\msyfam=\sevenmsy
  \scriptscriptfont\msyfam=\fivemsy

\def\hexnumber@#1{\ifcase#1 0\or1\or2\or3\or4\or5\or6\or7\or8\or9\or
	A\or B\or C\or D\or E\or F\fi }

\edef\msx@{\hexnumber@\msxfam}
\edef\msy@{\hexnumber@\msyfam}

\mathchardef\boxdot="2\msx@00
\mathchardef\boxplus="2\msx@01
\mathchardef\boxtimes="2\msx@02
\mathchardef\square="0\msx@03
\mathchardef\blacksquare="0\msx@04
\mathchardef\centerdot="2\msx@05
\mathchardef\lozenge="0\msx@06
\mathchardef\blacklozenge="0\msx@07
\mathchardef\circlearrowright="3\msx@08
\mathchardef\circlearrowleft="3\msx@09
\mathchardef\rightleftharpoons="3\msx@0A
\mathchardef\leftrightharpoons="3\msx@0B
\mathchardef\boxminus="2\msx@0C
\mathchardef\Vdash="3\msx@0D
\mathchardef\Vvdash="3\msx@0E
\mathchardef\vDash="3\msx@0F
\mathchardef\twoheadrightarrow="3\msx@10
\mathchardef\twoheadleftarrow="3\msx@11
\mathchardef\leftleftarrows="3\msx@12
\mathchardef\rightrightarrows="3\msx@13
\mathchardef\upuparrows="3\msx@14
\mathchardef\downdownarrows="3\msx@15
\mathchardef\upharpoonright="3\msx@16

\mathchardef\downharpoonright="3\msx@17
\mathchardef\upharpoonleft="3\msx@18
\mathchardef\downharpoonleft="3\msx@19
\mathchardef\rightarrowtail="3\msx@1A
\mathchardef\leftarrowtail="3\msx@1B
\mathchardef\leftrightarrows="3\msx@1C
\mathchardef\rightleftarrows="3\msx@1D
\mathchardef\Lsh="3\msx@1E
\mathchardef\Rsh="3\msx@1F
\mathchardef\rightsquigarrow="3\msx@20
\mathchardef\leftrightsquigarrow="3\msx@21
\mathchardef\looparrowleft="3\msx@22
\mathchardef\looparrowright="3\msx@23
\mathchardef\circeq="3\msx@24
\mathchardef\succsim="3\msx@25
\mathchardef\gtrsim="3\msx@26
\mathchardef\gtrapprox="3\msx@27
\mathchardef\multimap="3\msx@28
\mathchardef\therefore="3\msx@29
\mathchardef\because="3\msx@2A
\mathchardef\doteqdot="3\msx@2B

\mathchardef\triangleq="3\msx@2C
\mathchardef\precsim="3\msx@2D
\mathchardef\lesssim="3\msx@2E
\mathchardef\lessapprox="3\msx@2F
\mathchardef\eqslantless="3\msx@30
\mathchardef\eqslantgtr="3\msx@31
\mathchardef\curlyeqprec="3\msx@32
\mathchardef\curlyeqsucc="3\msx@33
\mathchardef\preccurlyeq="3\msx@34
\mathchardef\leqq="3\msx@35
\mathchardef\leqslant="3\msx@36
\mathchardef\lessgtr="3\msx@37
\mathchardef\backprime="0\msx@38
\mathchardef\risingdotseq="3\msx@3A
\mathchardef\fallingdotseq="3\msx@3B
\mathchardef\succcurlyeq="3\msx@3C
\mathchardef\geqq="3\msx@3D
\mathchardef\geqslant="3\msx@3E
\mathchardef\gtrless="3\msx@3F
\mathchardef\sqsubset="3\msx@40
\mathchardef\sqsupset="3\msx@41
\mathchardef\vartriangleright="3\msx@42
\mathchardef\vartriangleleft="3\msx@43
\mathchardef\trianglerighteq="3\msx@44
\mathchardef\trianglelefteq="3\msx@45
\mathchardef\bigstar="0\msx@46
\mathchardef\between="3\msx@47
\mathchardef\blacktriangledown="0\msx@48
\mathchardef\blacktriangleright="3\msx@49
\mathchardef\blacktriangleleft="3\msx@4A
\mathchardef\vartriangle="0\msx@4D
\mathchardef\blacktriangle="0\msx@4E
\mathchardef\triangledown="0\msx@4F
\mathchardef\eqcirc="3\msx@50
\mathchardef\lesseqgtr="3\msx@51
\mathchardef\gtreqless="3\msx@52
\mathchardef\lesseqqgtr="3\msx@53
\mathchardef\gtreqqless="3\msx@54
\mathchardef\Rrightarrow="3\msx@56
\mathchardef\Lleftarrow="3\msx@57
\mathchardef\veebar="2\msx@59
\mathchardef\barwedge="2\msx@5A
\mathchardef\doublebarwedge="2\msx@5B
\mathchardef\angle="0\msx@5C
\mathchardef\measuredangle="0\msx@5D
\mathchardef\sphericalangle="0\msx@5E
\mathchardef\varpropto="3\msx@5F
\mathchardef\smallsmile="3\msx@60
\mathchardef\smallfrown="3\msx@61
\mathchardef\Subset="3\msx@62
\mathchardef\Supset="3\msx@63
\mathchardef\Cup="2\msx@64

\mathchardef\Cap="2\msx@65

\mathchardef\curlywedge="2\msx@66
\mathchardef\curlyvee="2\msx@67
\mathchardef\leftthreetimes="2\msx@68
\mathchardef\rightthreetimes="2\msx@69
\mathchardef\subseteqq="3\msx@6A
\mathchardef\supseteqq="3\msx@6B
\mathchardef\bumpeq="3\msx@6C
\mathchardef\Bumpeq="3\msx@6D
\mathchardef\lll="3\msx@6E

\mathchardef\ggg="3\msx@6F

\mathchardef\circledS="0\msx@73
\mathchardef\pitchfork="3\msx@74
\mathchardef\dotplus="2\msx@75
\mathchardef\backsim="3\msx@76
\mathchardef\backsimeq="3\msx@77
\mathchardef\complement="0\msx@7B
\mathchardef\intercal="2\msx@7C
\mathchardef\circledcirc="2\msx@7D
\mathchardef\circledast="2\msx@7E
\mathchardef\circleddash="2\msx@7F
\def\ulcorner{\delimiter"4\msx@70\msx@70 }
\def\urcorner{\delimiter"5\msx@71\msx@71 }
\def\llcorner{\delimiter"4\msx@78\msx@78 }
\def\lrcorner{\delimiter"5\msx@79\msx@79 }
\def\yen{\mathhexbox\msx@55 }
\def\checkmark{\mathhexbox\msx@58 }
\def\circledR{\mathhexbox\msx@72 }
\def\maltese{\mathhexbox\msx@7A }
\mathchardef\lvertneqq="3\msy@00
\mathchardef\gvertneqq="3\msy@01
\mathchardef\nleq="3\msy@02
\mathchardef\ngeq="3\msy@03
\mathchardef\nless="3\msy@04
\mathchardef\ngtr="3\msy@05
\mathchardef\nprec="3\msy@06
\mathchardef\nsucc="3\msy@07
\mathchardef\lneqq="3\msy@08
\mathchardef\gneqq="3\msy@09
\mathchardef\nleqslant="3\msy@0A
\mathchardef\ngeqslant="3\msy@0B
\mathchardef\lneq="3\msy@0C
\mathchardef\gneq="3\msy@0D
\mathchardef\npreceq="3\msy@0E
\mathchardef\nsucceq="3\msy@0F
\mathchardef\precnsim="3\msy@10
\mathchardef\succnsim="3\msy@11
\mathchardef\lnsim="3\msy@12
\mathchardef\gnsim="3\msy@13
\mathchardef\nleqq="3\msy@14
\mathchardef\ngeqq="3\msy@15
\mathchardef\precneqq="3\msy@16
\mathchardef\succneqq="3\msy@17
\mathchardef\precnapprox="3\msy@18
\mathchardef\succnapprox="3\msy@19
\mathchardef\lnapprox="3\msy@1A
\mathchardef\gnapprox="3\msy@1B
\mathchardef\nsim="3\msy@1C
\mathchardef\ncong="3\msy@1D

\mathchardef\varsubsetneq="3\msy@20
\mathchardef\varsupsetneq="3\msy@21
\mathchardef\nsubseteqq="3\msy@22
\mathchardef\nsupseteqq="3\msy@23
\mathchardef\subsetneqq="3\msy@24
\mathchardef\supsetneqq="3\msy@25
\mathchardef\varsubsetneqq="3\msy@26
\mathchardef\varsupsetneqq="3\msy@27
\mathchardef\subsetneq="3\msy@28
\mathchardef\supsetneq="3\msy@29
\mathchardef\nsubseteq="3\msy@2A
\mathchardef\nsupseteq="3\msy@2B
\mathchardef\nparallel="3\msy@2C
\mathchardef\nmid="3\msy@2D
\mathchardef\nshortmid="3\msy@2E
\mathchardef\nshortparallel="3\msy@2F
\mathchardef\nvdash="3\msy@30
\mathchardef\nVdash="3\msy@31
\mathchardef\nvDash="3\msy@32
\mathchardef\nVDash="3\msy@33
\mathchardef\ntrianglerighteq="3\msy@34
\mathchardef\ntrianglelefteq="3\msy@35
\mathchardef\ntriangleleft="3\msy@36
\mathchardef\ntriangleright="3\msy@37
\mathchardef\nleftarrow="3\msy@38
\mathchardef\nrightarrow="3\msy@39
\mathchardef\nLeftarrow="3\msy@3A
\mathchardef\nRightarrow="3\msy@3B
\mathchardef\nLeftrightarrow="3\msy@3C
\mathchardef\nleftrightarrow="3\msy@3D
\mathchardef\divideontimes="2\msy@3E
\mathchardef\varnothing="0\msy@3F
\mathchardef\nexists="0\msy@40
\mathchardef\mho="0\msy@66
\mathchardef\eth="0\msy@67
\mathchardef\eqsim="3\msy@68
\mathchardef\beth="0\msy@69
\mathchardef\gimel="0\msy@6A
\mathchardef\daleth="0\msy@6B
\mathchardef\lessdot="3\msy@6C
\mathchardef\gtrdot="3\msy@6D
\mathchardef\ltimes="2\msy@6E
\mathchardef\rtimes="2\msy@6F
\mathchardef\shortmid="3\msy@70
\mathchardef\shortparallel="3\msy@71
\mathchardef\smallsetminus="2\msy@72
\mathchardef\thicksim="3\msy@73
\mathchardef\thickapprox="3\msy@74
\mathchardef\approxeq="3\msy@75
\mathchardef\succapprox="3\msy@76
\mathchardef\precapprox="3\msy@77
\mathchardef\curvearrowleft="3\msy@78
\mathchardef\curvearrowright="3\msy@79
\mathchardef\digamma="0\msy@7A
\mathchardef\varkappa="0\msy@7B
\mathchardef\hslash="0\msy@7D
\mathchardef\hbar="0\msy@7E
\mathchardef\backepsilon="3\msy@7F

\catcode`\@=12

\def\eqalign#1{\null\,\vcenter{\openup\jot \mathsurround=0pt
  \ialign{\strut\hfil$\displaystyle{##}$&$
        \displaystyle{{}##}$\hfil \crcr#1\crcr}}\,}

\def\epsilon{\varepsilon}
\def\sgn{\mathop{\rm sgn}}
\def\twosum#1#2{\sum_{\scriptstyle #1 \atop \scriptstyle #2}}

\def\phi{\varphi}

\def\C{{\cal{C}}}
\def\F{{\cal{F}}}
\def\G{{\cal{G}}}
\def\H{{\cal{H}}}
\def\O#1{{\cal O }(#1)}
\def\U{{\cal{U}}}
\def\I{{\cal{I}}}
\def\J{{\cal{J}}}
\def\K{{\cal{K}}}
\def\semi{\mathbin{\circledS}}

\def\R{\BB{R}}
\def\Z{\BB{Z}}
\def\i{{\MM{i}}}
\def\j{{\MM{j}}}
\def\l{{\MM{l}}}
\def\balpha{{\MM{\alpha}}}
\def\si{{\Mm{i}}}
\def\sj{{\Mm{j}}}
\def\sl{{\Mm{l}}}
\def\sbalpha{{\Mm{\alpha}}}

\begin{document}

\title[SPATIAL DISCRETIZATION OF PDEs WITH INTEGRALS]{Spatial 
discretization of partial differential equations with integrals}

\author{Robert I. McLachlan}
\address{\hskip-\parindent
        Mathematics\\ 
        Massey University\\
        Private Bag 11--222\\
        Palmerston North\\
		New Zealand
		}
\email{R.McLachlan@massey.ac.nz}

\thanks{Research at MSRI is supported in part by NSF grant DMS-9701755.}

\begin{abstract}
We consider the problem of constructing spatial finite
difference approximations on a fixed, arbitrary
grid, which have analogues
of any number of integrals of the partial differential
equation and of some of its symmetries. 
A basis for the
space of of such difference operators is constructed;
most cases of interest involve a single such basis element.
(The ``Arakawa'' Jacobian is such an element.)
We show how the topology of the grid affects the complexity
of the operators. 
\end{abstract}

\maketitle

\section{Conservative discretization}

{\narrower \it ``Numerical methods for nonlinear conservation laws are among
the great success stories of modern numerical analysis.'' \rm
(Iserles \cite{iserles})

}

\bigskip
Conservative discretizations of partial differential equations have
been explored for a long time. What does ``conservative'' mean?
An early definition is due to Lax and Wendroff \cite{la-we}, who
considered the class on PDEs with one spatial dimension,
\begin{equation}
\label{conspde}
u_t + \partial_x(f(u)) = 0 
\end{equation}
and called discretizations of the form
\begin{equation}
\label{laxwendroff}
{u^{n+1}_i - u^n_i \over \Delta t} =
{H(u^n_{i+j},\dots,u^n_{i-j+1}) - H(u^n_{i+j-1},\dots,u^n_{i-j})\over
\Delta x}
\end{equation}
conservative. See \cite{iserles} for an introduction to such
methods. More generally, the formulation (\ref{conspde})
is called conservative, and the expanded form
$$ u_t + f'(u)u_x $$
nonconservative, with these terms carrying over to the corresponding
discrete forms. The PDE (\ref{conspde}) reflects, amongst other
things, conservation of the integral of $u$ (e.g. total mass,
momentum, etc.), and the conservative discretization
(\ref{laxwendroff}) preserves a discrete analog:
$\sum_i u^{n+1}_i = \sum_i u^n$. The full consequences
for the discrete scheme of the form (\ref{laxwendroff})
remain unclear.

More recently the term has been applied to PDEs that can be 
written purely in terms of intrinsic differential operators such as 
div, grad, and curl. A conservative spatial discretization is
then one which preserves discrete analogues of these operators' integral
identities (e.g., Stokes's theorem). In many cases these obey maximum
principles and have robust stability properties in difficult situations
such as rough grids and discontinuous coefficients \cite{shashkov}.

Schemes have also been developed for particular equations that inherit
conserved quantities approximating those of the PDE. An early and
famous example is the Arakawa Jacobian \cite{arakawa}, a discretization
of $v_x w_y - v_y w_x$ which, when applied to the two-dimensional
Euler fluid equations, provides two conservation laws corresponding
to energy and enstrophy, both quadratic functions. It is widely used
in computational meteorology. There are many energy-conserving
schemes for particular PDEs: Fei and V\'asquez \cite{fe-va} for
the sine-Gordon equation; Glassey \cite{glassey} for the 
Zakharov equations; 
Glassey and Schaeffer \cite{gl-sc} for a nonlinear wave equation. 
The original presentations of all these are somewhat ad-hoc, the
proof of conservation relying on a telescoping sum.

The Arakawa Jacobian has the extremely nice property that
it can be applied to systems (in two space dimensions,
with two variables) with {\it any} two integrals, 
not just energy and enstrophy. It was further explained
and generalized to arbitrary grids by Salmon and Talley
in \cite{sa-ta}. It is this systematic approach that
we generalize in this paper to equations with any number
of integrals, space dimensions, and variables. Our formulation
includes all integral-preserving discretizations.

Having integrals of course reduces the evolution to a smaller
space, which, when their level sets are compact, gives
the method a form or nonlinear stability. Often, more is true:
Preservation of a discrete form of $\int u\,dx$ by the
Lax-Wendroff form (\ref{laxwendroff}) leads to correct
shock speeds, and preservation of energy and enstrophy 
by the Arakawa Jacobian prevents energy cascading to small 
length scales \cite{arakawa}. 

The Euler equations, the sine-Gordon equation and so on are examples
of {\it Hamiltonian} PDEs, which suggests that one should look for
semi- or fully-discrete forms which preserve not only a discrete energy 
but also a discrete Hamiltonian (symplectic) structure. 
For systems with canonical Hamiltonian structure, this possibility
was explored in \cite{mclachlan}.  But
even before the importance of Hamiltonian PDEs was widely recognized,
for which a watershed event was perhaps the 1983 conference \cite{sandiego},
it had been pointed out by Morrison \cite{phil} that spatial
discretizations of non-canonical Hamiltonian PDEs will not normally
be Hamiltonian. One example apart, the curious `sine bracket' Hamiltonian 
discretization of the Jacobian \cite{zeitlin}, this
is a difficult and essentially unsolved problem. 
Probably the right
generalization of `Hamiltonian' has not yet been found. 

We are thus
reluctantly led to consider only energy-conserving discretizations.
Or perhaps we should not be reluctant: Simo et al. 
\cite{si-ta-wo} have argued
and presented detailed evidence from elastodynamics that conserving
energy leads to excellent nonlinear stability properties that 
preserving symplectic structure does not. (Essentially because
symplectic schemes can only brake the fast modes, whereas energy-conserving
schemes can also damp them.)

In Hamiltonian systems, energy is normally seen as playing a distinguished
role. Yet there may be other conserved quantities just as important
for the long-time dynamics. Some of them, the `Casimir' integrals,
can be due to the Hamiltonian structure itself. Non-Hamiltonian
systems can also have conserved quantities. Even in the ODE example
of the free rigid body,
the relationship between
schemes preserving energy and/or momentum and/or symplectic structure
is quite complicated \cite{le-si}.

In this paper we go some way towards uniting these different 
integrals and different points of view.
Our goal is to develop a methodology for building spatial discretizations
that preserve discrete analogues of any given set of integrals.
It should be systematic, all-inclusive, and reproduce known schemes.
We do this in a formulation in which the integrals appear explicitly;
the integrals themselves can then be discretized in any way. The 
basic ``finite difference molecule'' is now a completely skew-symmetric
tensor. Symmetry plays a fundamental role, and we will see how the
skew-symmetry of this tensor interacts with other desired symmetries
of the scheme (e.g., translational and rotational invariance) in 
a nontrivial way.

\section{Hamiltonian and other PDEs with integrals}

We consider PDEs with independent spatial variables $x\in\R^d$ and 
dependent variables $u(x)\in\R^m$. We loosely call $m$ the ``number
of variables.'' The relevant class of sufficiently smooth
real-valued functionals of $u(x)$ will be denoted $\U$.
A Hamiltonian PDE is specified by a Hamiltonian
$\H\in\U$ and a Poisson bracket $\{\, , \, \}:\U\times \U\to \U$:
\begin{equation}
\label{eom}
\dot u = \{u,\H\},
\end{equation}
where the Poisson bracket is bilinear, skew-symmetric
\begin{equation}
\{\F,\G\} = -\{\G,\F\},
\end{equation}
and obeys the Jacobi identity
\begin{equation}
\{\F,\{\G,\H\}\} + \{\G,\{\H,\F\}\} + \{\H,\{\F,\G\}\} = 0
\end{equation}
and the Leibniz rule
\begin{equation}
\{\F,\G\H\} = \{\F,\G\}\H + \{\F,\H\}\G
\end{equation}
for all $\F$, $\G$, $\H\in\U$.
In fact, these axioms imply the existence of a Hamiltonian (or `Poisson')
operator $\J$, such that the Poisson bracket takes the form
\begin{equation}
\{\F,\G\} = \int {\delta \F \over \delta u} \J {\delta\G\over \delta u}\,dx
\end{equation}
and (\ref{eom}) becomes
\begin{equation}
\dot u = \J {\delta \H \over \delta u}
\end{equation}
A particularly important example is the conservation law (\ref{conspde}),
which takes this form with $x$, $u\in R$ and
\begin{equation}
\label{conspde2}
\H = \int F(u(x))\,dx,\quad \J = \partial_x,\quad F' = f.
\end{equation}

The system (\ref{eom}) has functional $\I\in\U$ as an integral
if $\dot \I = \{\I,\H\}$=0. Some integrals $\C$ are distinguished in
that $\{C,\F\}=0$ for all $\F\in \U$; they are called Casimirs. 
The operator $\J=\partial_x$ has a single Casimir, $\C = \int u\,dx$,
because $\J(\delta \C/\delta u) = \partial_x 1 = 0$.

The peculiarly Hamiltonian character of these PDEs is due to
the Jacobi identity satisfied by the Poisson bracket. If, as
argued previously, we discard this identity, what class of systems
result? Energy is still conserved, because of the skew-symmetry
of the bracket. Systems (\ref{eom}) may still have integrals
and the operator $\J$ may still have Casimirs. Do such systems
retain any other special properties? The answer is no: {\it all} systems
with an integral $\H$ can be written in the form (\ref{eom}) for
some choice of the bracket (or equivalently, for some choice
of the skew-adjoint operator $\J$). So to study systems with
an integral, and their discretizations, we may without
loss of generality assume the form (\ref{eom}).

This is most easily seen in the finite dimensional case, as shown
recently by Quispel and Capel \cite{qu-ca}. Take
$u\in \R^n$ as coordinates on phase space. Poisson ($\equiv$
non-canonical Hamiltonian) systems 
$$ \dot u = \{u,H\} = J(u)\nabla H(u) $$
have integral $H$.
But suppose an arbitrary system $\dot u = f(u)$ has integral $H$.
Let $z=\nabla H$ and
\begin{equation}
\label{Jij}
J_{ij} = {f_i z_j - f_j z_i\over \sum_k z_k^2}
\end{equation}
Then $J$ is skew-symmetric, and $J\nabla H=f$ as required. 
This $J$ is singular at critical points of $H$, but in \cite{mc-qu-ro} it
is shown that if $f$ and $H$ are smooth, and the critical points
of $H$ are nondegenerate, then there is a smooth matrix $J$ such
that $f=J\nabla H$.

This idea extends easily to systems with any number of integrals:
The system of ODEs $\dot u = f(x)$
has integrals $I^1,\dots,I^p$ if and only if there exists a
totally skew-symmetric $(p+1)$-tensor $K$ such that for all
$x$ where the vectors $\nabla I^i$ are linearly independent,
\begin{equation}
\label{Kijk}
 f_i = K_{ijk\dots} {\partial I^1\over \partial x_j}
{\partial I^2\over\partial x_k}\dots. 
\end{equation}
For, suppose $K$ exists. Then $\dot I^j = f\cdot\nabla I^j
=0$, so each $I^j$ is an integral. Conversely, suppose $f$ has integrals
$I^j$. Then (using exterior algebra, see \cite{darling})
$$K={f\wedge\nabla I^1\dots\wedge\nabla I^p \over 
\det(\nabla I^i\cdot\nabla I^j)}
$$
satisfies (\ref{Kijk}).
$K$ is determined uniquely only in the case $n=p+1$; 
see \cite{mc-qu-ro} for further details.
We write the inner product (\ref{Kijk}) as
$$
f = K(\nabla I^1,\nabla I^2,\dots). 
$$

What about Casimirs? Suppose that instead of contracting $K$ against
all the integrals, as in (\ref{Kijk}), we contract against just one,
say $I^1$. Then $\widetilde K = K(\nabla I^1)$ is a skew $p$-tensor which has 
$I^1$ as a Casimir, in the sense that $\widetilde K(\nabla I^1) \equiv 0$ and
any differential equation formed from this $\widetilde K$ (as in (\ref{Kijk}))
will have $I^1$ as an integral. But there are many different $K$'s, all
generating the same system $\dot u = f$, that do not have $I^1$ as a Casimir.
Thus the distinction between the Hamiltonian, other integrals, and Casimirs,
that was present for Hamiltonian systems, is lost now. We are free
to move between different representations of $f$ as needed.

(There is another importance difference. If $J\in \R^{n\times n}$ 
satisfies the Jacobi identity and has locally constant rank $m$, then
$J$ automatically has $n-m$ Casimirs \cite{olver}. This need not be true if the
Jacobi identity does not hold: there may not be $n-m$ functions whose
gradients span $J$'s nullspace. This is another reason for constructing
$J$'s as above that automatically have the required Casimirs.)

One can use the tensor $K$ to define a $(p+1)$-bracket,
$$ \{F_1,\dots,F_{p+1}\} = K(\nabla F_1,\dots,\nabla F_{p+1}) $$
which is multilinear, a derivation in each argument, and completely
antisymmetric. Such brackets have been revived in modern times
by Nambu \cite{nambu}, who, amongst other things, introduced the
3-bracket on $\R^3$ given by $K_{ijk} = \epsilon_{ijk}$ (the
alternating tensor). This gives systems of the form
$$ \dot u = K(\nabla I^1,\nabla I^2) = \nabla I^1 \times \nabla I^2. $$
In particular, the free rigid body takes this form with
$I^1={1\over 2}|u|^2$ being total angular momentum and
$I^2={1\over 2}\sum u_i^2/A_i$ being the kinetic energy, where
the $A_i$ are the body's moments of inertia.
Contracting against $\nabla I^1$ gives the standard, Lie-Poisson form
of the equations,
$$ \dot u_i = J_{ik}(u) {\partial I^2\over u_k} = \epsilon_{ijk} u_j 
{\partial I^2 \over \partial u_k}.$$
However, later studies \cite{gautheron,takhtajan}, attempting to build
a true generalization of Hamiltonian mechanics from such brackets,
have found that not all constant $K$'s satisfy the required ``fundamental
identity'' (the analogue of the Jacobi identity); its solutions
all lead to systems with $n-1$ integrals. So it is not clear that
interesting dynamics as well as interesting algebraic structure
will arise in this way.

The situation for PDEs is formally the same. For example, for a PDE
$\dot u = f(u)$
with one integral $\I$, one can define the operator $\cal J$
analogously to (\ref{Jij}) by
$$ \J(x,x')v(x') = {\int \left(f(x)z(x') - f(x')z(x)\right)v(x')\,dx' \over
\int {\delta \I\over\delta u}(x')z(x')\,dx' }$$
where $z(x)$ is any smooth function, and
$$ f(x) = \int \J(x,x'){\delta I\over\delta u}(x')\,dx'. $$
The extension to multiple integrals $\{\I^j\}$ is similar.
However, questions of convergence of the integrals arise, 
and the nonuniqueness
situation is much worse: it is not clear how to construct
{\it local} operators, for example. However, as our goal is to construct
finite-dimensional finite difference operators, the representation
(\ref{Kijk}) is sufficient.

Looking back at Eq. (\ref{conspde}), we see that it encompasses two
important conservation laws expressed in two different ways.
The fact that $\partial_x (\delta u/\delta u) = \partial_x 1 = 0$
means that the Casimir $\int u\,dx$ is an integral, and the
conservative scheme (\ref{laxwendroff}) maintains a discrete
analog of this. 
In this example, this property is relatively easy to preserve
under discretization: a system has $\sum_i u_i$ as an integral
if and only if it can be written in the form $\dot u_i = J_{ij} F_j$,
where $J$ is not necessarily skew, but $\sum_i J_{ij}=0$ for all $j$.
Without loss of generality we
can take $J$ to be in the form (\ref{laxwendroff}), with just two
nonzero diagonals. The form of the $F_j$ chosen in (\ref{laxwendroff})
is necessary for translation invariance.

Secondly, if $f$ is a variational derivative,
$f=\delta \H/\delta u$ say, then the skew-adjointness of 
$\J = \partial_x$ means that $\H$ is an integral.
To preserve this property under discretization means taking (again,
without loss of generality) $\dot u_i = J_{ij}{\partial H\over \partial
u_j}$, where $J$ is skew symmetric. Note that such a $J$ need
not {\it a priori} have $\sum u_i$ as a Casimir, and, similarly,
the nonsymmetric $J$ used in (\ref{laxwendroff}) does not preserve
any discrete $H$.
Thus, the two expressions of conservation laws are in fact independent.

In this paper we generalize the second form. The first form is
deceptively simple in this example, because the Casimir is so simple.
It is not clear how to modify (\ref{laxwendroff}) to incorporate
different Casimirs. In the second form the integral appears
explicitly and, once $J$ is found, {\it any} quantity can be conserved.

Before continuing, we mention one trivial but complete solution to
the whole problem. Why not contract $K$ against {\it all} the
integrals and have simply $\dot u_i = f_i$, with all $I^j$ 
being integrals of $f$? That is, $f\cdot\nabla I^j=0$ for all $j$.
This case is already included in the above formulation, with
$f$ regarded as a 
``skew 1-tensor.'' The equations $f\cdot\nabla I^j=0$ are linear
and can be solved in many ways, for example, by starting with
an arbitrary $f$ and projecting to the subspace $\{f\cdot \nabla I^j=0\:
j=1,\dots,p\}$. One objection is that this solution is
so general that it is not clear how to proceed in any particular case.
For example, to modify $f$ as little as possible one might choose orthogonal
projection, but this will couples all of the $f_i$.
By incorporating more of the known structure of the problem
we can work more systematically.

\section{Method of discretization}

We wish to construct discretizations of the form (\ref{Kijk}).
There are two ways to proceed. One could take a particular PDE,
write it in the form $\dot u = \K(\delta I^1,\dots)$,
and discretize this skew-adjoint operator $\K$, preserving
skew symmetry. This is difficult, if only because such formulations
of PDEs are new and have not been widely developed yet.
Instead, we study systems of the
form (\ref{Kijk}) in their own right, constructing
elementary tensors $K_{ijk\dots}$, and
seeing what PDEs they can be used to approximate. That is, we establish
(in a sense defined below)
$$ K(v^1,\dots,v^p) \dots = \K(v^1,\dots,v^p) + \O{h^r}.$$
We call $K$ a {\it finite difference tensor} and $K(v^1,\dots,v^p)$
a {\it finite difference}.
Then, the integrals $\I^i$ can be discretized in any way, say by
$$ I^i(u) = \I^i(u) + \O{h^r}$$
(this amounts to a numerical quadrature)
giving the conservative system of ODEs
\begin{equation}
\label{Kbracket}
\dot u = K(\nabla I^1,\dots,\nabla I^p).
\end{equation}
Since this form includes all systems with integrals $I^i$, we can be
confident of not missing any in our construction.

We start with an elementary example illustrating how
easy it is to break skew symmetry. With one integral in 
one space dimension, we are seeking an antisymmetric matrix $K$.
On a constant-spaced grid, central differences have
\begin{equation}
\label{central}
K = {1\over 2h}\left(
\begin{array}{cccc}
\ddots & & & \\
-1 & 0 & 1 & \\
& -1 & 0 & 1 \\
& & & \ddots \\
\end{array}
\right),
\end{equation}
which is antisymmetric. On a non-constant-spaced grid,
if we let $p(x)$ be the quadratic interpolating
$(x_0,v_0)$, $(x_1,v_1)$, and $(x_2,v_2)$ and use
the estimate $v_1'= p'(x_1)$, the associated matrix
$K$ is not antisymmetric---it even has a nonzero diagonal.
Taking its antisymmetric part is not a good idea, as we
have no idea what operator $K^{\rm T}$ approximates.
Indeed, it is not immediately clear what bandwidth is
required to achieve order 2, say, with an antisymmetric matrix.
The element $K_{i,i-1}$, which is relevant to $\dot u_i$, must
also contribute to $\dot u_{i-1}$.

Nonconstant operators also pose problems. Let $\K = u\partial_x
+ \partial_x u$. Central differences are not skew,
but it is not obvious that the skew matrix 
\begin{equation}
\label{nonskew}
K = {1\over h}\left(\begin{array}{cccc}
\ddots & & & \\
-u_{i-1} & 0 & u_i & \\
& -u_i & 0 & u_{i+1} \\
& & & \ddots \\
\end{array}\right),
\end{equation}
is a discretization of $\K$, or how to increase its order from 1.

Below we develop some requirements on the tensors $K$, and construct
all the elementary ones, for various numbers of integrals, space
dimensions, and grids.

\section{Definitions \& theory}

The fundamental objects are the {\it grid} $L$, the
{\it index set} $M$, the {\it symmetry group} $G$, 
and the skew tensors $K\in \Lambda^{p+1}(\R^{M})$, which
we now define.

Let $L$ be a set of indices of grid points. To each index $i\in L$
there is a physical point $x_i\in\R^d$. Let $\{1,\dots,m\}$ be
the set of indices of the dependent variables, so that the full,
discrete state space is indexed by the {\it index set}
$$ M = L\times\{1,\dots,m\}.$$
A {\it grid function} is a real function on $M$; 
for example, the system
state is given by the grid function $u:M\to\R$. Its value at
point $i=(l,\alpha)\in M$ is written $u_i=u_{(l,\alpha)}$.
(That is, we are assembling all the unknowns into a big ``column
vector.'')
When $m=1$, we drop the second subscript entirely.

For simplicity, we only consider the 
interpretation of this function in which 
$u_{(l,\alpha)}\approx u_\alpha(x_l)$.
(Staggered grids and $u_i$ representing other functionals of $u(x)$ do not 
affect our main line of argument.)

The $p$ functions, and their corresponding grid functions,
which are to be inserted in (\ref{Kbracket}) are denoted $v^1,\dots,v^p$.
However, when $p=1$ we denote it $v$, and when $p=2$
we denote them $v$ and $w$, to reduce the number of indices.

A discretization of a PDE is thus a vector field on $\R^M$,
and a discretization of an operator $\K$ is a skew $(p+1)$-tensor
$K_{i_0\dots i_p} = K_\si= K_{(\sl,\sbalpha)}
\in \R$, $i_j\in M$; i.e.,
$K\in \Lambda^{p+1}(\R^M)$. Since we are trying to {\it construct}
such tensors, intermediate steps will also involve nonskew
tensors, i.e., real functions on $M^{p+1}$.

$K$ {\it approximates $\K$ to order $r$} if 
$$ K_{i_0\dots i_p}v^1(x_{i_1})\dots v^p(x_{i_p}) =
\K(v_1,\dots,v^p)(x_{i_0})+\O{h^r}$$
for all smooth $v^1,\dots,v^p$.
We sometimes 
drastically abbreviate this to $Kv = \K v + \O{h^r}$.
We also abbreviate $v^1_{i_1}\dots v^p_{i_p}$ to $v_\si$.

Let the grid $L$ have a nonnegative distance 
function $|j-k|$. This
extends to $M$  by $|(j,\alpha)-(k,\beta)| := |j-k|$.
The {\it bandwidth} of $K$ is the smallest $c$ such that
$K_\si=0$ for all $\i$ such that $|i_j - i_k|>c$.

Two examples are 
the Euclidean distance $|x_i-x_j|$, giving the ``Euclidean
bandwidth,'' and 
the minimum number of edges traversed going from $i$ and $j$,
where the grid points have been connected to form a graph,
giving the ``graph bandwidth.''
For example, $K$ in Eq. (\ref{nonskew}) has graph bandwidth 1.
$K$ is {\it local} if it has a finite bandwidth even on 
infinite grids.

The grid $L$ may be structured, like
a square or triangular lattice, or unstructured. Let $G$ be
a symmetry group acting on $M=L\times\{1,\dots,m\}$. 
We usually consider only {\it spatial symmetries}, those which are
the identity in their second slot, i.e. $\pi_2 g(l,\alpha)=\alpha$.
These merely rearrange grid functions on $L$. 

Furthermore, $G$ is usually a subgroup of the
symmetry group of the continuous physical space.
For example, suppose this space is the plane. Many PDEs of
physical interest are invariant under the group $E(2)$ 
of Euclidean motions of the plane
(the
semi-direct product of rotations, reflections, and translations,
$O(2)\semi\R^2$). 
A discrete version of such a PDE can inherit
some of this invariance if $L$ has a subgroup of $E(2)$ as a symmetry
group. Examples are the square lattice, which has $D_4\semi\Z^2$
as a symmetry group (8 rotations and reflections, plus discrete translations)
and the equilateral triangular lattice, which has $D_3\semi\Z^2$.
For the Euclidean group of the sphere there are no such natural lattices,
and the dislocations that occur, e.g., 
when triangulating an icosahedral grid,
are known to destabilize numerical methods and create artifacts in the 
solutions.

Some nonspatial symmetries can also be included. For example,
for a PDE involving $f(v_{2x}-v_{1y})$ we might include the map
$(x,y,v_1,v_2)\mapsto (y,x,v_2,v_1)$ in $G$.

Order of accuracy can also sometimes be expressed as a symmetry. 
One way to ensure second-order accuracy is for the expansion
of the discretization error in powers of the spatial grid size $h$
to have only odd or only even terms present. 
This is equivalent to being invariant
under the operation $h\mapsto -h$, or $x\mapsto -x$. This can only
apply if $x\mapsto-x$ is a symmetry of the lattice itself, which it
is for square and triangular lattices.

To include this possibility we equip each element of $G$ with a 
sign, $\sgn(g)=\pm 1$, such that $G$ is homomorphic to $\Z_2$.
The map corresponding to $h\mapsto -h$ would then have sign
$1$ ($-1$) when the operator has an even (odd) number of derivatives.

The action of $G$ extends to an action on $\Lambda^{p+1}(\R^M)$ by
$$g K_{i_0\dots i_p} := K_{g(i_0)\dots g(i_p)}$$
which we write as
$$(g K)_\si =  g(K_\si) = K_{g(\si)}.$$
A tensor $K$ is {\it $G$-invariant} if $g K = (\sgn g) K$ for all
$g\in G$.

Thus we have the following requirements on the finite difference operator $K$:
\begin{itemize}
\item $K$ should be completely skew-symmetric;
\item $K$ should be $G$-invariant;
\item $K$ should be as simple as possible;
\item $K$ should approximate the desired continuous operator to
the desired order;
\item $K$ should be local.
\end{itemize}
However, these requirements conflict with each other.

One way to construct operators such as $K$  
with the required symmetry properties
is to sum over the symmetry group. Given any 
tensor $K$, 
\begin{equation}
\label{sum1}
\sum_{\sigma\in S_{p+1},g\in G} \sgn(\sigma)\sgn(g) 
K_{\sigma(g(\si))} 
\end{equation}
is completely skew-symmetric and $G$-invariant.
This suggests two ways to construct symmetric $K$'s:
\begin{itemize}
\item Start with a $K$ which approximates the desired continuous operator,
and symmetrize it;
\item Start with a very simple $K$, such as a basis element for the
space of $(p+1)$-tensors, symmetrize it, and see what continuous operator it
approximates.
\end{itemize}
A major drawback of the first strategy is that we have no control over
what the symmetrized $K$ approximates. 

The second strategy builds a ``library'' of all such
difference operators, from which linear combinations can 
be taken as desired. However, the form (\ref{sum1}) is not convenient
for writing down these operators in the usual way, which
requires the coefficients of each $v^j_k$ appearing
in the resulting grid function at a particular point
$i_0$. That is, we want to know $K_{i_0 i_1 \dots i_p}$ for
a particular $i_0$.  We derive finite differences in this
form in three stages: firstly, for $m=1$ variable;
secondly, for $m\ge 1$ variables with no unknowns at the
same grid point coupled; thirdly, the general case, $m\ge 1$
variables with arbitrary coupling.

\subsection*{Case 1. $m=1$ variable}
When $m=1$ we drop the second component
of the index $i=(l,\alpha)$.
We can take $i_0=0$.
Fix a multi-index $\i\in M^{p+1}$ where $i_0=0$ 
and start with the elementary tensor defined by $K_\si = 1$, 
$K_\sj=0$
for all $\j \ne \i$. We assume that the indices in $\i$
are distinct, for otherwise skew-symmetrizing $K$ would lead to
the zero tensor. Skew-symmetrizing $K$ gives a tensor of bandwidth 
$\max_{j,k}|i_j-i_k|$. It is
$$
H_\sl = \sum_{\rho\in S_{p+1},h\in G}
\sgn(\rho) \sgn(h)  K_{h(\rho(\sl))}
$$
so the vector field at the point $0$ is $\sum_\sl H_\sl v_\sl$
where $l_0=0$. Since $K$ is a discrete delta function,
there is only one nonzero term in this sum, i.e., 
\begin{equation}
\label{Hsl}
\sum_\sl H_\sl v_\sl = \sum_{\rho,h,\sl}
\sgn(\rho)\sgn(h) v_\sl,
\end{equation}
where the sum is taken over all $\rho$, $h$, and $\l$ 
such that $\rho(h(\l))=\i$ and $l_0=0$. Therefore,
$\l = \rho^{-1}(h^{-1}(\i))$. Let $g=h^{-1}$, $\sigma=\rho^{-1}$
so that Eq. (\ref{Hsl}) becomes
\begin{equation}
\label{sv}
\sum_{\sigma,g}\sgn(\sigma)\sgn(g) v_{\sigma(g(\si))}
\end{equation}
where the sum is over all $\sigma$ and $g$ such that
$\sigma(g(\i))_0=0$.

Since no two indices in $\i$ are the same, 
for each such $g$, let $\sigma$ be such that
$\sigma(g(\i))_0=0$; then the remaining $\sigma$'s 
that satisfy this equation 
lie in $S_p$, the
permutations of the last $p$ indices. The sum over $S_p$
can be evaluated to give a determinant, giving the
vector field at point 0,
\begin{equation}
\label{the_f}
F(\i) := \sum_{g\in G_\si}\sgn (\sigma)\sgn (g)\det V_{\sigma(g(\si))}
\end{equation}
where 
\begin{equation}
\label{Gi}
G_\si = \{ g\in G: \exists \sigma\in S_{p+1}\hbox{\rm\ such that\ }
\sigma(g(\i))_0 = 0\},
\end{equation}
and the $p\times p$ matrix $V$ has $(j,k)$ entry
\begin{equation}
\label{V}
(V_\sl)_{jk} = v^j_{l_k},\quad 1\le j,k\le p.
\end{equation}
(The first subscript of $l_0=0$
and does not appear in the matrix.)

We introduce a graphical notation for formulas such as Eq. 
(\ref{the_f}). 
An arrow connecting gridpoints
$i_1,\dots,i_p$ will indicate a term $\det V_\si$.
For $p>1$, the sign factors can be incorporated by applying 
a permutation of
sign $\sgn(\sigma)\sgn(g)$ to the $i_j$ (for $p=2$ and $p=3$
we merely change the direction of the arrow if the required sign
is $-1$, equivalent to writing the columns of $V$ in reverse order)
or by choosing $\sigma$ in (\ref{Gi}) so that $\sgn(\sigma)\sgn(g)=1$.
The reader is encouraged
to refer immediately to Fig. (\ref{fig1}(a)) and its
associated finite difference Eq. (\ref{p2d1}) to see how
easy this is.

Thus, constructing skew finite differences amounts to choosing an
initial arrow, finding the group $G_\si$, and finding the image
of the initial arrow under $G_\si$. 

\subsection*{Case 2. $m\ge 1$ variables, distinct points
coupled}
Let $\i=(\l,\balpha)$ is the chosen basis element,
where $l_0=0$. ``Distinct points coupled'' means that
$l_j\ne l_k$ for all $j$, $k$. 

The single element $\i$ will contribute to the vector field at 
$(0,\alpha_j)$ for all $j$; 
therefore, we will construct, not a single basis vector field,
but the family of vector fields spanned by $K_{(\sl,\sbalpha)}$
{\it for all} $\balpha$. That is, we allow coupling of all 
components right from the start. We introduce the rank $p+1$,
dimension $m$ tensor $T_\sbalpha
\in \R^{m^{p+1}}$
and start with $K_{\si,\sbalpha}=T_\sbalpha$, $K_{\sj,\sbalpha}=0$
for $\j\ne \i$.

Passing from 
Eq. (\ref{sv}) to Eq. (\ref{the_f}) only required that
the $l_j$ be distinct; therefore the
symmetrized vector field at the point 0
(corresponding to Eq. (\ref{the_f})
in the single variable case) is
\begin{equation}
\twosum{g\in G_\si}{\alpha_0,\dots,\alpha_p}\sgn (\sigma)\sgn (g)
T_\sbalpha\det V_{\sigma(g(\si))}
{\partial \over \partial u_{\sigma(g(\sbalpha))_0}}
\end{equation}
We want to find the $\alpha_0$-component of this vector field.
To do this we relabel the
dummy indices $\balpha$ by applying
$(\sigma g)^{-1}$ {\it in the second slot only} to get
the vector field at point 0 in component $\alpha_0$,
\begin{equation}
\label{the_f2}
F(\i) := \twosum{g\in G_\si}{\alpha_1,\dots,\alpha_p}
\sgn (\sigma)\sgn (g)
T_{g^{-1}(\sigma^{-1}(\sbalpha))}\det V_{(\sigma(g(\sl)),\sbalpha)}
\end{equation}

Notice that in (\ref{the_f2}), each determinant involves
the same components of the $v^i$. Also, if $g$ is a spatial symmetry,
then $g^{-1}(\sigma^{-1}(\balpha))=\sigma^{-1}(\balpha)$.

The diagram notation extends easily to (\ref{the_f2}). 
To the arrow $\sigma(g(\l))$ we attach the label
$g^{-1}(\sigma^{-1}(\balpha))$ indicating the $T$-tensor
attached to that determinant.

\subsection*{Case 3. $m\ge 1$ variables, arbitrary coupling.}

Equality of some of the $l_j$ affects the sum over permutations
in Eq. (\ref{sv}). 
Let $\j=h(\rho(\l))$ where $j_0=0$. Let $n(\j)$ be the
number of 0's in $\j$. Then the subgroup of $S_{p+1}$
leaving $\sigma(\j)_0=0$ is not $S_p$ as it was before.
It is convenient to have a sum of determinants of
all $p\times p$ matrices, so we express this subgroup as
the product of $S_p$ and the flips $(0k)$, $k=1,\dots,n(\j)-1$.
This subgroups of $S_{p+1}$ 
contains all permutations of the first $n$
elements; summing over these merely skew-symmetrizes $T$.
We could have imposed this in the first place for simplicity.
Since there will usually be $g_j$ (a translation, say)
with $g_j(l_j)=0$, this is true for any set of equal
elements of $\l$.

Let $\Sigma_l=\{\sigma: \sigma(\l)=\l\}$ be the symmetry
group of $\l$. To sum up, we
\begin{equation}
\label{Tsym}
\hbox{\rm take $T$ to skew-symmetric under $\Sigma_\sl$.}
\end{equation}

With this assumption, each flip $(0k)$ gives
an equal contribution and we can evaluate the sum
over permutations to give the symmetrized vector
field at the point 0
\begin{equation}
\label{the_f3}
F(\i) := \twosum{g\in G_\si}{\alpha_1,\dots,\alpha_p}
\sgn (\sigma)\sgn (g) n(g(\l))
T_{g^{-1}(\sigma^{-1}(\sbalpha))}\det V_{(\sigma(g(\sl)),\sbalpha)}
\end{equation}

In the diagrams, to the arrow $\sigma(g(\l))$ we attach
the weight $n(g(\l))$.

To summarize, the final finite difference evaluated at
the point $(0,\alpha_0)$, is given by Equations (\ref{the_f3}) 
(\ref{Gi}), and (\ref{V}).  Eq. (\ref{the_f3}) specializes
to Eq. (\ref{the_f2}) when all elements of $\l$ are distinct,
and specializes further to (\ref{the_f}) when $m=1$.
In practice, from a diagram
one writes down the finite difference directly from its
diagram.

We develop these diagrams and study the resulting differences 
for different numbers of
integrals and dimensions of phase space, and different
symmetry groups $G$.

\section{Examples}

\subsection*{Case 1. $p=1$ integral, $d=1$ space dimension.}
By scaling it is sufficient to consider $\i=(0,1)$.
The bandwidth is 1.
The group $G_\si$ has the single element $i\mapsto i-1$. The
permutation which brings $0$ to the front is $(0,-1)\mapsto (-1,0)$,
with sign $-1$. Thus we get the standard central difference
$$F(0,1) = v_1 - v_{-1}.$$
On a grid with constant spacing $h$, 
$$ v_1 - v_{-1} = 2 h v_x + \O{h^3}$$.

All other examples are related to this one:
\begin{enumerate}
\item by Richardson extrapolation,
$$ 8 F(0,1) - F(0,2)  = 12 h v_x + \O{h^5}.$$
\item eliminating the leading order term(s) gives finite differences 
approximating
higher-order differential operators:
$$ F(0,2) - 2 F(0,1) = 2 h^3 v_{xxx} + \O{h^5}.$$
\item taking a linear combination of these basis elements gives stencils that
approximate
other first-order differential operators. This is equivalent to multiplying 
by a symmetric tensor $s_{ij}$. To get a smooth continuous
limit we can take, e.g., 
$s_{ij} = q(x_i,x_j,u_i,u_j)$, where $q$ is symmetric in its first and
second pairs of arguments. If $K$ is the tensor corresponding to 
the vector field $F(0,1)$,
\begin{equation}
\label{gen1d}
s_{ij}K_{ij} v_j = 2h(s\partial_x + 
\partial_x s)v + \O{h^3},
\end{equation}
where $s := q(x,x,u,u)$ (no sum on $i$).
\end{enumerate}
We consider this last example in more detail, since it gives the class
of all skew tridiagonal (i.e., bandwidth 1) finite differences.

Firstly, suppose we have the non-constant operator $\J = u\partial_x
+ \partial_x u$. Eq. (\ref{gen1d}) discretizes this if $s(x,u) = u$.
The only constraint is the the tensor $s_{ij}$ must be symmetric,
to maintain the overall skew-symmetry of the finite difference.
For example, we can take $s_{ij} = (u_i+u_j)/2$ (only elements with
$|i-j|=1$ are actually used). This gives the finite difference
tensor
\begin{equation}
\label{vp}
K = {1\over 2}\left(\begin{array}{ccccc}
\ddots & & & &\\
-u_0-u_1 & 0 & u_1 + u_2 & & \\
& -u_1 - u_2 & 0 & u_2 + u_3 & \\
& & -u_2 - u_3 & 0 & u_3 + u_4 \\
& & & & \ddots \\
\end{array}\right),
\end{equation}
showing how the skew-symmetry is maintained.
\footnote{Interestingly, the choice $s_{ij} = \sqrt{u_i u_j}$ 
actually gives a {\it Poisson} $K$. This is because it is the
image of the standard central difference under a
change of variables which sends $\partial_x$ 
into $u \partial_x + \partial_x u$.}

Secondly,
suppose we wish to difference on an irregular grid,
where the data are known at the points $c(x_i)= c(ih)$---we know
$v(c(x_i))$. Then we
want an approximation of $v' = v_c$. From the chain rule, this is 
equal to $v_x/c_x$. We cannot get
this by applying (\ref{gen1d}) to $v$, since the terms
in $v$ only cannot cancel. 

We apply (\ref{gen1d}) 
to a function $w(x,v(c(x))$. This  gives the equation 
$$ (s\partial_x + \partial_x s)w(x,v(c(x))) =
(\sqrt{2s}\partial_x\sqrt{2s})w(x,v(c(x)) = {v_x(c(x))\over c'(x)}$$
with solution
$$s = {1\over 2{c'}^2},\quad w = c' v.$$

For the discretization, any $w_j$ and any symmetric $s_{ij}$ 
with this continuous
limit can be taken. If $c'$ is known analytically, we can use the
midpoints of the intervals to get
$$s_{0,1} = {1\over2(c_{1\over2}')^2},\quad 
s_{0,-1}={-1\over 2(c_{-{1\over2}}')^2},$$
$$ v'_0 = {1\over 4h}\left({c'_1 v_1 \over {c_{1\over2}'}^2 }- 
{c'_{-1} v_{-1} \over {c_{-{1\over2}}'}^2} \right)
+\O{h^2}$$
---that is, a second order, 2-point, anti-symmetric discretization of the
derivative on smooth grids.
With constant spacing, $c(x)=x$ and it reduces to the standard
central difference.

If $c'$ is not known, we can approximate it symmetrically with
$$f_{0i} = {1\over2}\left({c_i-c_0\over h}\right)^{-2},$$
$$w_i = {(v_{i}+v_{i+1})(c_{i+1}-c_i) + (v_i+v_{i-1})(c_i-c_{i-1})\over4h},$$
which has bandwidth 2. It is not consistent on rough grids, however;
the error is $\O{h^2 c''''}=\O{h^{-1}}$ near a discontinuity of $c'$.

With $m>1$ components and bandwidth 1, the two possibilities
are $\l = (0,0)$ and $\l=(0,1)$. In the first case, the tensor
$T$ must be skew symmetric in its only two slots, but there
are no group symmetries. In the second case, $G_\sl$
has two elements, the identity and left translation. Under
left translation,
$(0,1)\mapsto (-1,0)$; applying $\sigma=(01)$ of sign $-1$ maps
$(-1,0)\mapsto(0,-1)$. In 
the last step we apply the $\sigma^{-1}$ to
the indices of $T$. Combining both possibilities gives the
finite difference
$$ -T^{\rm T} v_{-1} + J v_0 + T v_1 = (J+T-T^{\rm T})v + 
h (T+T^{\rm T})v_x +\O{h^2},$$
where $J=-J^{\rm T}$.
Further imposing the symmetry $i\to -i$, of sign $-1$, is equivalent
to taking $T=T^{\rm T}$; then the finite difference is second order.
As above, nonconstant operators are approximated by taking
$T_l=T(x_l, x_{l+1}, u_l,u_{l+1})$, symmetric
in its first and second pairs of arguments, and
$J_l = J(x_{l-1},x_l,x_{l+1},u_{l-1},u_l,u_{l+1})$,
symmetric under $(13)$ and $(46)$.

 \begin{figure}[tb]
 \begin{center}
 \leavevmode

\setlength{\unitlength}{4144sp}%
\begingroup\makeatletter\ifx\SetFigFont\undefined%
\gdef\SetFigFont#1#2#3#4#5{%
  \reset@font\fontsize{#1}{#2pt}%
  \fontfamily{#3}\fontseries{#4}\fontshape{#5}%
  \selectfont}%
\fi\endgroup%
\begin{picture}(4861,1225)(1074,-2011)
\thinlines
\put(4065,-1109){\circle*{70}}
\put(4508,-1101){\circle*{70}}
\put(5408,-1101){\circle*{70}}
\put(5858,-1109){\circle*{70}}
\put(5528,-1109){\vector( 1, 0){225}}
\put(5273,-1116){\vector(-1, 0){653}}
\put(4178,-1109){\vector( 1, 0){210}}
\put(5385,-899){\makebox(0,0)[lb]{\smash{\SetFigFont{9}{10.8}{\rmdefault}{\mddefault}{\updefault}
1
}}}
\put(5828,-891){\makebox(0,0)[lb]{\smash{\SetFigFont{9}{10.8}{\rmdefault}{\mddefault}{\updefault}
2
}}}
\put(4935,-900){\makebox(0,0)[lb]{\smash{\SetFigFont{9}{10.8}{\rmdefault}{\mddefault}{\updefault}
0
}}}
\put(3961,-913){\makebox(0,0)[lb]{\smash{\SetFigFont{9}{10.8}{\rmdefault}{\mddefault}{\updefault}
-2
}}}
\put(4418,-914){\makebox(0,0)[lb]{\smash{\SetFigFont{9}{10.8}{\rmdefault}{\mddefault}{\updefault}
-1
}}}
\put(4951,-1366){\makebox(0,0)[b]{\smash{\SetFigFont{9}{10.8}{\rmdefault}{\mddefault}{\updefault}
cab
}}}
\put(5626,-1366){\makebox(0,0)[b]{\smash{\SetFigFont{9}{10.8}{\rmdefault}{\mddefault}{\updefault}
abc
}}}
\put(4951,-1591){\makebox(0,0)[b]{\smash{\SetFigFont{9}{10.8}{\rmdefault}{\mddefault}{\updefault}
cab
}}}
\put(4276,-1591){\makebox(0,0)[b]{\smash{\SetFigFont{9}{10.8}{\rmdefault}{\mddefault}{\updefault}
abc
}}}
\put(4276,-1366){\makebox(0,0)[b]{\smash{\SetFigFont{9}{10.8}{\rmdefault}{\mddefault}{\updefault}
bca
}}}
\put(5626,-1546){\makebox(0,0)[b]{\smash{\SetFigFont{9}{10.8}{\rmdefault}{\mddefault}{\updefault}
bca
}}}
\put(1178,-1141){\circle*{70}}
\put(1621,-1133){\circle*{70}}
\put(2521,-1133){\circle*{70}}
\put(2971,-1141){\circle*{70}}
\put(2641,-1141){\vector( 1, 0){225}}
\put(2386,-1148){\vector(-1, 0){653}}
\put(1291,-1141){\vector( 1, 0){210}}
\put(2498,-931){\makebox(0,0)[lb]{\smash{\SetFigFont{9}{10.8}{\rmdefault}{\mddefault}{\updefault}
1
}}}
\put(2941,-923){\makebox(0,0)[lb]{\smash{\SetFigFont{9}{10.8}{\rmdefault}{\mddefault}{\updefault}
2
}}}
\put(2048,-932){\makebox(0,0)[lb]{\smash{\SetFigFont{9}{10.8}{\rmdefault}{\mddefault}{\updefault}
0
}}}
\put(1074,-945){\makebox(0,0)[lb]{\smash{\SetFigFont{9}{10.8}{\rmdefault}{\mddefault}{\updefault}
-2
}}}
\put(1531,-946){\makebox(0,0)[lb]{\smash{\SetFigFont{9}{10.8}{\rmdefault}{\mddefault}{\updefault}
-1
}}}
\put(1966,-2011){\makebox(0,0)[lb]{\smash{\SetFigFont{10}{12.0}{\rmdefault}{\mddefault}{\updefault}
(a)
}}}
\put(4906,-1996){\makebox(0,0)[lb]{\smash{\SetFigFont{10}{12.0}{\rmdefault}{\mddefault}{\updefault}
(b)
}}}
\end{picture}
 \caption{Case 2, Two integrals, one space dimension. (a) One variable;
 (b) $m$ variables. See Eq. (\ref{p2d1}) for the finite difference
 interpretation of (a).}
 \label{fig1}
 \end{center}
 \end{figure}
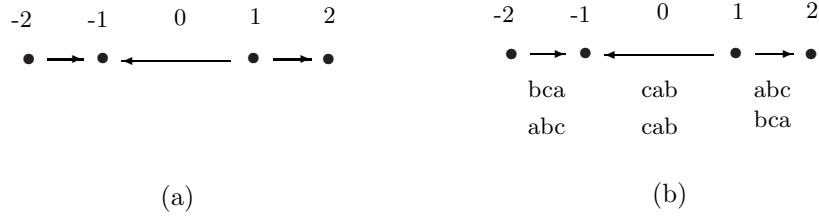

\subsection*{Case 2. $p=2$ integrals, $d=1$ space dimension}

With $m=1$ component, the simplest tensor has base index
$\i=(0,1,2)$. This will lead to bandwidth 2.
We start with the arrow $1\to 2$ (see Figure \ref{fig1}). 
Translating left by 1 and
rotating indices right (an even permutation) (i.e.,
performing $(0,1,2)\mapsto (-1,0,1)\mapsto (0,1,-1)$)
gives the arrow $1\to -1$. Repeating gives the
arrow $-2\to-1$. The resulting diagram is already
symmetric under $h\to -h$, so we do not need to add
this operation.
The diagram in Fig. \ref{fig1}(a) corresponds to the finite
difference
\begin{equation}
\label{p2d1}
\eqalign{
F(0,1,2) &= 
\det V_{1,2} + \det V_{1,-1} + \det V_{-2,-1} \cr
& = (v_1 w_2-v_2 w_1) + (v_{1}w_{-1}-v_{-1}w_{1}) +
(v_{-2}w_{-1}-v_{-1}w_{-2})
}
\end{equation}
Expanding in Taylor series, this is
$$ h^3(3(v'w''-w'v'') + 2(vw'''-wv''')) + \O{h^5}$$

With $m>1$ component there is essentially one finite difference
each with bandwidth 0, 1, and 2. 
We write $T_{\alpha_1,\alpha_2,\alpha_3}
=T_{abc}$.

With $\l=(0,0,0)$ (bandwidth 0), Eq. (\ref{Tsym}) 
says we must have $T_{abc}$
completely skew-symmetric. One might not call this a ``difference,''
since it only acts on $v_0$ and $w_0$.

 \begin{figure}[tb]
 \begin{center}
 \leavevmode
 
\setlength{\unitlength}{4144sp}%
\begingroup\makeatletter\ifx\SetFigFont\undefined%
\gdef\SetFigFont#1#2#3#4#5{%
  \reset@font\fontsize{#1}{#2pt}%
  \fontfamily{#3}\fontseries{#4}\fontshape{#5}%
  \selectfont}%
\fi\endgroup%
\begin{picture}(4250,827)(1341,-406)
\thinlines
\put(1671, 89){\oval(170,122)[tr]}
\put(1671, 59){\oval(182,182)[tl]}
\put(1671, 59){\oval(182,182)[bl]}
\put(1671, 29){\oval(170,122)[br]}
\put(1756, 29){\vector( 0, 1){0}}
\put(5131, 29){\vector( 0, 1){0}}
\put(5216, 29){\oval(170,122)[bl]}
\put(5216, 59){\oval(182,182)[br]}
\put(5216, 59){\oval(182,182)[tr]}
\put(5216, 89){\oval(170,122)[tl]}
\put(1891, 74){\circle*{70}}
\put(2791, 74){\circle*{70}}
\put(2341, 74){\circle*{70}}
\put(4996, 74){\circle*{70}}
\put(4096, 74){\circle*{70}}
\put(4546, 74){\circle*{70}}
\put(2431, 74){\vector( 1, 0){225}}
\put(4456, 74){\vector(-1, 0){225}}
\put(2775,308){\makebox(0,0)[lb]{\smash{\SetFigFont{9}{10.8}{\rmdefault}{\mddefault}{\updefault}
1
}}}
\put(2325,307){\makebox(0,0)[lb]{\smash{\SetFigFont{9}{10.8}{\rmdefault}{\mddefault}{\updefault}
0
}}}
\put(4996,299){\makebox(0,0)[lb]{\smash{\SetFigFont{9}{10.8}{\rmdefault}{\mddefault}{\updefault}
1
}}}
\put(4546,299){\makebox(0,0)[lb]{\smash{\SetFigFont{9}{10.8}{\rmdefault}{\mddefault}{\updefault}
0
}}}
\put(4051,299){\makebox(0,0)[b]{\smash{\SetFigFont{9}{10.8}{\rmdefault}{\mddefault}{\updefault}
--1
}}}
\put(1890,316){\makebox(0,0)[b]{\smash{\SetFigFont{9}{10.8}{\rmdefault}{\mddefault}{\updefault}
--1
}}}
\put(1441,-83){\makebox(0,0)[b]{\smash{\SetFigFont{9}{10.8}{\rmdefault}{\mddefault}{\updefault}
bca
}}}
\put(2558,-188){\makebox(0,0)[b]{\smash{\SetFigFont{9}{10.8}{\rmdefault}{\mddefault}{\updefault}
abc(2)
}}}
\put(4344,-151){\makebox(0,0)[b]{\smash{\SetFigFont{9}{10.8}{\rmdefault}{\mddefault}{\updefault}
abc(2)
}}}
\put(5491,-76){\makebox(0,0)[b]{\smash{\SetFigFont{9}{10.8}{\rmdefault}{\mddefault}{\updefault}
bca
}}}
\put(4494,-391){\makebox(0,0)[lb]{\smash{\SetFigFont{10}{12.0}{\rmdefault}{\mddefault}{\updefault}
(b)
}}}
\put(2251,-406){\makebox(0,0)[lb]{\smash{\SetFigFont{10}{12.0}{\rmdefault}{\mddefault}{\updefault}
(a)
}}}
\end{picture}
 \caption{Case 2, Two integrals, one space dimension, $m$ variables,
 $\l=(0,0,1)$.}
 \label{fig2}
 \end{center}
 \end{figure}
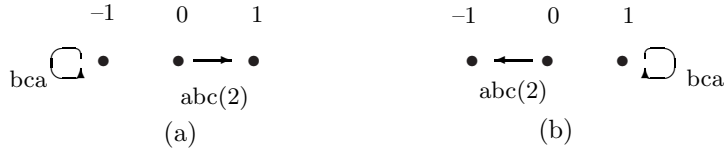

With $\l=(0,0,1)$ (bandwidth 1), Eq. (\ref{Tsym})
says we have $T_{abc}=-T_{bac}$.
The weight of $\l$ is $n(\l)=2$, since it has 2 zeros.
$G_\sl$ has two elements, the identity and a left translation.
Applying the  left translation followed by 
a shift-right permutation $\sigma$ (of sign 1), 
$(0,0,1)\mapsto(-1,-1,0)\mapsto(0,-1,-1)$, giving 
the arrow $-1\to-1$ with label $\sigma^{-1}(abc)=bca$. 
$(0,-1,-1)$ has one zero, so the weight of this arrow is 1.
Together we get two arrows, with diagram 
Fig. \ref{fig2}(a) and finite difference
$$ \eqalign{& 2 T_{abc} (v_{0,b}w_{1,c}-v_{1,c}w_{0,b})
+ T_{bca}(v_{-1,b}w_{-1,c}-v_{-1,c}w_{-1,b}) \cr
& = (T_{abc}+T_{bca}+T_{cab})v_b w_c + 
h (T_{abc}-T_{bca})v_b w'_c +
h (T_{cab} - T_{bca})v'_b w_c + {\cal O}(h^2).
}
$$
(The initial, skew-symmetric term could be removed by
a term $F(0,0,0)$.)

The resulting tensor is not invariant under $g:i\mapsto -i$.
Its image under $g$ is shown in Fig. \ref{fig2}(b).
These two diagrams can be added or subtracted to get a tensor
that is $g$-invariant with sign $1$ or $-1$, as desired.  

With $\l=(0,1,2)$ (bandwidth 2), Eq. (\ref{Tsym}) says
that $T$ is arbitrary.
The group $G_\sl$ has three elements: the identity, and a shift
left by 1 or 2.
Apply the two translations gives the diagram Fig. \ref{fig1}(b). 
However,
unlike this case with $m=1$,
this is not invariant under $i\mapsto -i$, i.e., it
does not give a second-order finite difference. Applying
this symmetry gives the second row of labels in Fig. \ref{fig3}(b).
(For example, under $(0,1,2)\mapsto
(0,-1,-2)$ the arrow
$1\to 2(abc)$ maps to the arrow
$-1\to -2(abc)$ with sign $-1$, or
$-2\to -1(abc)$ with sign 1. Adding these
makes $T_{abc}=T_{bca}$,
i.e., we can take $T$ to be symmetric under even permutations.

\subsection*{Case 3. $p=1$ integral, $d=2$ space dimensions}
With $p$ free indices in $K$ we can only couple unknowns which span
a $p$-dimensional subspace of $\R^d$. This is equivalent to
the case $d=p$. For example, on a square grid in $\R^2$,
$F((0,0),(0,1))= h v_y + {\cal O}(h^3)$.

Thus, to get fundamentally new finite difference tensors, we need
$p\ge d$.

 \begin{figure}[tb]
 \begin{center}
 \leavevmode
 
\setlength{\unitlength}{4144sp}%
\begingroup\makeatletter\ifx\SetFigFont\undefined%
\gdef\SetFigFont#1#2#3#4#5{%
  \reset@font\fontsize{#1}{#2pt}%
  \fontfamily{#3}\fontseries{#4}\fontshape{#5}%
  \selectfont}%
\fi\endgroup%
\begin{picture}(3377,3323)(1343,-3098)
\thinlines
\put(1870,-82){\circle*{70}}
\put(2320,-75){\circle*{70}}
\put(1420,-82){\circle*{70}}
\put(1420,-532){\circle*{70}}
\put(2312,-532){\circle*{70}}
\put(2320,-982){\circle*{70}}
\put(1870,-975){\circle*{70}}
\put(1870,-525){\circle*{70}}
\put(1420,-967){\circle*{70}}
\put(1750, 90){\makebox(0,0)[lb]{\smash{\SetFigFont{9}{10.8}{\rmdefault}{\mddefault}{\updefault}
(0,1)
}}}
\put(2470,-563){\makebox(0,0)[lb]{\smash{\SetFigFont{9}{10.8}{\rmdefault}{\mddefault}{\updefault}
(1,0)
}}}
\put(4120,-82){\circle*{70}}
\put(4570,-75){\circle*{70}}
\put(3670,-82){\circle*{70}}
\put(3670,-532){\circle*{70}}
\put(4562,-532){\circle*{70}}
\put(4570,-982){\circle*{70}}
\put(4120,-975){\circle*{70}}
\put(4120,-525){\circle*{70}}
\put(3670,-967){\circle*{70}}
\put(4000, 90){\makebox(0,0)[lb]{\smash{\SetFigFont{9}{10.8}{\rmdefault}{\mddefault}{\updefault}
(0,1)
}}}
\put(4720,-563){\makebox(0,0)[lb]{\smash{\SetFigFont{9}{10.8}{\rmdefault}{\mddefault}{\updefault}
(1,0)
}}}
\put(4478,-450){\vector(-1, 1){277}}
\put(3676,-181){\vector( 0,-1){255}}
\put(4208,-983){\vector( 1, 0){248}}
\put(1870,-1882){\circle*{70}}
\put(2320,-1875){\circle*{70}}
\put(1420,-1882){\circle*{70}}
\put(1420,-2332){\circle*{70}}
\put(2312,-2332){\circle*{70}}
\put(2320,-2782){\circle*{70}}
\put(1870,-2775){\circle*{70}}
\put(1870,-2325){\circle*{70}}
\put(1420,-2767){\circle*{70}}
\put(4120,-1882){\circle*{70}}
\put(4570,-1875){\circle*{70}}
\put(3670,-1882){\circle*{70}}
\put(3670,-2332){\circle*{70}}
\put(4562,-2332){\circle*{70}}
\put(4570,-2782){\circle*{70}}
\put(4120,-2775){\circle*{70}}
\put(4120,-2325){\circle*{70}}
\put(3670,-2767){\circle*{70}}
\put(2217,-432){\vector(-1, 1){277}}
\put(4186,-2701){\vector( 1, 1){270}}
\put(4036,-1958){\vector(-1,-1){285}}
\put(4471,-1891){\vector(-1, 0){240}}
\put(3668,-2416){\vector( 0,-1){232}}
\put(3766,-2783){\vector( 1, 0){240}}
\put(1478,-2416){\vector( 1,-1){293}}
\put(2311,-2686){\vector( 0, 1){248}}
\put(1763,-1891){\vector(-1, 0){232}}
\put(2228,-2250){\vector(-1, 1){277}}
\put(1958,-2783){\vector( 1, 0){248}}
\put(1418,-1981){\vector( 0,-1){255}}
\put(3728,-2416){\vector( 1,-1){293}}
\put(4561,-2686){\vector( 0, 1){248}}
\put(4013,-1891){\vector(-1, 0){232}}
\put(4478,-2250){\vector(-1, 1){277}}
\put(4208,-2783){\vector( 1, 0){248}}
\put(3668,-1981){\vector( 0,-1){255}}
\put(4568,-2228){\vector( 0, 1){240}}
\put(1801,-1298){\makebox(0,0)[lb]{\smash{\SetFigFont{9}{10.8}{\rmdefault}{\mddefault}{\updefault}
(a)
}}}
\put(4066,-1261){\makebox(0,0)[lb]{\smash{\SetFigFont{10}{12.0}{\rmdefault}{\mddefault}{\updefault}
(b)
}}}
\put(1748,-3098){\makebox(0,0)[lb]{\smash{\SetFigFont{10}{12.0}{\rmdefault}{\mddefault}{\updefault}
(c)
}}}
\put(4118,-3091){\makebox(0,0)[lb]{\smash{\SetFigFont{10}{12.0}{\rmdefault}{\mddefault}{\updefault}
(d)
}}}
\put(1750,-1710){\makebox(0,0)[lb]{\smash{\SetFigFont{9}{10.8}{\rmdefault}{\mddefault}{\updefault}
(0,1)
}}}
\put(2470,-2363){\makebox(0,0)[lb]{\smash{\SetFigFont{9}{10.8}{\rmdefault}{\mddefault}{\updefault}
(1,0)
}}}
\put(4000,-1710){\makebox(0,0)[lb]{\smash{\SetFigFont{9}{10.8}{\rmdefault}{\mddefault}{\updefault}
(0,1)
}}}
\put(4720,-2363){\makebox(0,0)[lb]{\smash{\SetFigFont{9}{10.8}{\rmdefault}{\mddefault}{\updefault}
(1,0)
}}}
\end{picture}
 \caption{Case 4, Step-by-step construction of the Arakawa Jacobian.}
 \label{fig3}
 \end{center}
 \end{figure}
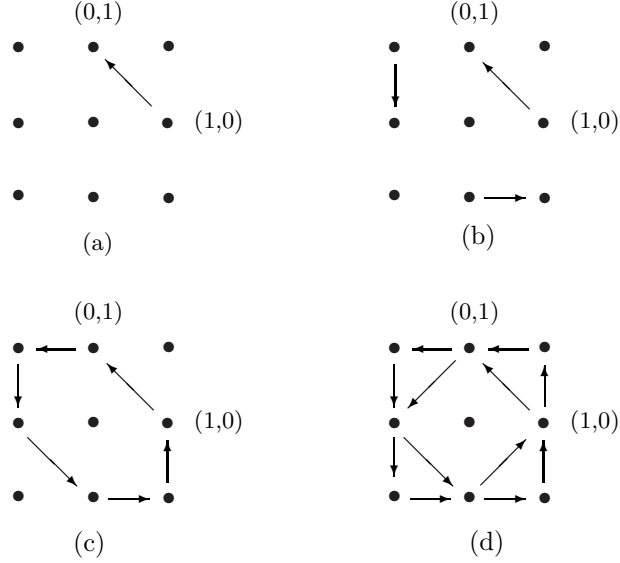

 \begin{figure}[tb]
 \begin{center}
 \leavevmode

\setlength{\unitlength}{4144sp}%
\begingroup\makeatletter\ifx\SetFigFont\undefined%
\gdef\SetFigFont#1#2#3#4#5{%
  \reset@font\fontsize{#1}{#2pt}%
  \fontfamily{#3}\fontseries{#4}\fontshape{#5}%
  \selectfont}%
\fi\endgroup%
\begin{picture}(5792,1500)(863,-1366)
\thinlines
\put(901,-526){\circle*{70}}
\put(1351,-518){\circle*{70}}
\put(1801,-518){\circle*{70}}
\put(1133,-136){\circle*{70}}
\put(1583,-128){\circle*{70}}
\put(1576,-908){\circle*{70}}
\put(1126,-908){\circle*{70}}
\put(2446,-519){\circle*{70}}
\put(2896,-511){\circle*{70}}
\put(3346,-511){\circle*{70}}
\put(2678,-129){\circle*{70}}
\put(3128,-121){\circle*{70}}
\put(3121,-901){\circle*{70}}
\put(2671,-901){\circle*{70}}
\put(4231,-181){\vector(-2,-3){174.462}}
\put(4320,-893){\vector( 1, 0){315}}
\put(4898,-428){\vector(-2, 3){164.769}}
\put(4736,-833){\vector( 2, 3){174.462}}
\put(4647,-121){\vector(-1, 0){315}}
\put(4069,-586){\vector( 2,-3){164.769}}
\put(4036,-511){\circle*{70}}
\put(4486,-503){\circle*{70}}
\put(4936,-503){\circle*{70}}
\put(4268,-121){\circle*{70}}
\put(4718,-113){\circle*{70}}
\put(4711,-893){\circle*{70}}
\put(4261,-893){\circle*{70}}
\put(5806,-174){\vector(-2,-3){174.462}}
\put(5895,-886){\vector( 1, 0){315}}
\put(6473,-421){\vector(-2, 3){164.769}}
\put(6311,-826){\vector( 2, 3){174.462}}
\put(6222,-114){\vector(-1, 0){315}}
\put(5644,-579){\vector( 2,-3){164.769}}
\put(5611,-504){\circle*{70}}
\put(6061,-496){\circle*{70}}
\put(6511,-496){\circle*{70}}
\put(5843,-114){\circle*{70}}
\put(6293,-106){\circle*{70}}
\put(6286,-886){\circle*{70}}
\put(5836,-886){\circle*{70}}
\put(1756,-443){\vector(-2, 3){164.769}}
\put(3301,-428){\vector(-2, 3){164.769}}
\put(2648,-188){\vector(-2,-3){174.462}}
\put(2723,-901){\vector( 1, 0){315}}
\put(1351,-1366){\makebox(0,0)[b]{\smash{\SetFigFont{10}{12.0}{\rmdefault}{\mddefault}{\updefault}
(a)
}}}
\put(2911,-1358){\makebox(0,0)[b]{\smash{\SetFigFont{10}{12.0}{\rmdefault}{\mddefault}{\updefault}
(b)
}}}
\put(4493,-1366){\makebox(0,0)[b]{\smash{\SetFigFont{10}{12.0}{\rmdefault}{\mddefault}{\updefault}
(c)
}}}
\put(1703, -1){\makebox(0,0)[b]{\smash{\SetFigFont{9}{10.8}{\rmdefault}{\mddefault}{\updefault}
(0,1)
}}}
\put(1966,-398){\makebox(0,0)[b]{\smash{\SetFigFont{9}{10.8}{\rmdefault}{\mddefault}{\updefault}
(1,0)
}}}
\put(6053, -8){\makebox(0,0)[b]{\smash{\SetFigFont{9}{10.8}{\rmdefault}{\mddefault}{\updefault}
bca
}}}
\put(5574,-263){\makebox(0,0)[b]{\smash{\SetFigFont{9}{10.8}{\rmdefault}{\mddefault}{\updefault}
cab
}}}
\put(5536,-781){\makebox(0,0)[b]{\smash{\SetFigFont{9}{10.8}{\rmdefault}{\mddefault}{\updefault}
abc
}}}
\put(6046,-1081){\makebox(0,0)[b]{\smash{\SetFigFont{9}{10.8}{\rmdefault}{\mddefault}{\updefault}
bca
}}}
\put(6526,-788){\makebox(0,0)[b]{\smash{\SetFigFont{9}{10.8}{\rmdefault}{\mddefault}{\updefault}
cab
}}}
\put(6555,-286){\makebox(0,0)[b]{\smash{\SetFigFont{9}{10.8}{\rmdefault}{\mddefault}{\updefault}
abc
}}}
\put(6068,-1358){\makebox(0,0)[b]{\smash{\SetFigFont{10}{12.0}{\rmdefault}{\mddefault}{\updefault}
(d)
}}}
\end{picture}
 \caption{Case 4, The Arakawa Jacobian on a triangular grid.}
 \label{fig4}
 \end{center}
 \end{figure}
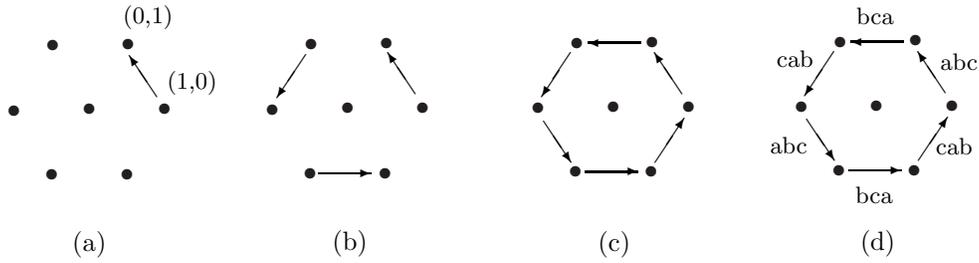

\subsection*{Case 4. $p=2$ integrals, $d=2$ space dimensions}
Consider $m=1$ and a square grid. The simplest index set
$\i$ we can take is $\i=((0,0),(1,0),(0,1))$, as shown
in Fig. \ref{fig3}(a). Unfortunately, this has lattice bandwidth 2
and Euclidean bandwidth $\sqrt{2}$, an unavoidable property
of the lattice.

Applying the two translations gives the diagram \ref{fig3}(b):
the simplest translation-invariant skew tensor. It gives
$h(v_x w_y - v_y w_x) + \O{h^2}$. That is, it is an
``Arakawa''-type Jacobian. 

With even operators on this grid, reflections (necessary
for second order accuracy) and rotations coincide,
which reduces the complexity of the generated finite
difference.
Applying them gives
the diagram \ref{fig3}(c), a second-order Jacobian. 
Finally, applying the rotations by $\pi/2$ gives
the diagram \ref{fig3}(d), a Jacobian with the full symmetry
group $D_4\semi \Z^2$.
As can be shown by expanding the entire finite difference,
{\it Fig. \ref{fig3}(d) is the Arakawa Jacobian} (first
derived in \cite{arakawa}.)

We could have stopped at Fig. \ref{fig3}(c); its anisotropy
may be irrelevant for some problems, and its complexity
is half that of \ref{fig2}(d)---12 terms instead of 24.

Consider the same problem on a regular triangular grid.
Now $\i=((0,0),(1,0),(0,1))$ (Fig. \ref{fig4}(a))
will give a graph bandwidth
of 1, not 2. Applying the two translations gives
Fig. \ref{fig4}(b), and reflections 
Fig. \ref{fig4}(c), which
{\it already} has the full symmetry of the grid. Thus
Arakawa-type Jacobians are naturally suited to triangular
grids. (Notice that Figs. \ref{fig4}(c) and \ref{fig3}(c)
are essentially the same.)

There are two points to learn from this:
\begin{enumerate}
\item 
With $p$ integrals, grids with $p+1$ mutual nearest neighbours
around a cell will give tensors of bandwidth 1. This
is only possible in dimension $d\ge p$.
\item
On some grids, the (optional) spatial symmetries coincide
with some of the (required) skew symmetries and/or reflection
symmetry (needed for second order accuracy).
\end{enumerate}
Cases 5 and 6 illustrate these points.

(An $m$-variable analogue of the Arakawa Jacobian is
shown in Fig. \ref{fig4}(c). It approximates a complicated
second order operator, but if $T_{abc}$ is symmetric
under even permutations, it is
$3\sqrt{3}h^2 T_{abc} \J(v_b,w_c)+{\cal O}(h^4)$.)

 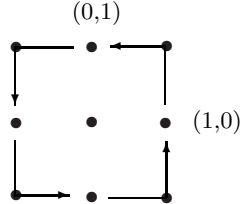
\begin{figure}[tb]
 \begin{center}
 \leavevmode

\setlength{\unitlength}{4144sp}%
\begingroup\makeatletter\ifx\SetFigFont\undefined%
\gdef\SetFigFont#1#2#3#4#5{%
  \reset@font\fontsize{#1}{#2pt}%
  \fontfamily{#3}\fontseries{#4}\fontshape{#5}%
  \selectfont}%
\fi\endgroup%
\begin{picture}(1127,1284)(893,-744)
\thinlines
\put(1420,233){\circle*{70}}
\put(1870,240){\circle*{70}}
\put(970,233){\circle*{70}}
\put(970,-217){\circle*{70}}
\put(1862,-217){\circle*{70}}
\put(1870,-667){\circle*{70}}
\put(1420,-660){\circle*{70}}
\put(1420,-210){\circle*{70}}
\put(970,-652){\circle*{70}}
\put(1300,405){\makebox(0,0)[lb]{\smash{\SetFigFont{9}{10.8}{\rmdefault}{\mddefault}{\updefault}
(0,1)
}}}
\put(2020,-248){\makebox(0,0)[lb]{\smash{\SetFigFont{9}{10.8}{\rmdefault}{\mddefault}{\updefault}
(1,0)
}}}
\put(1861,-106){\line( 0, 1){345}}
\put(1861,239){\vector(-1, 0){330}}
\put(1313,232){\line(-1, 0){352}}
\put(961,232){\vector( 0,-1){345}}
\put(961,-323){\line( 0,-1){330}}
\put(961,-653){\vector( 1, 0){322}}
\put(1523,-668){\line( 1, 0){345}}
\put(1868,-668){\vector( 0, 1){337}}
\end{picture}
 \caption{Case 5, 3 integrals in 2 dimensions}
 \label{fig5}
 \end{center}
 \end{figure}

\subsection*{Cases 5 \& 6: $p=3$ integrals in 2 \& 3 dimensions.}

The above observation suggests that in two dimensions,
the square grid, with 4 vertices
around each cell, is better suited to the case of three rather
than two 
integrals. With $\i=((0,0),(1,0),(1,1),(0,1))$, applying
the 3 translations only gives a tensor which is
$D_4$-symmetric (Fig. \ref{fig5}). It equals
$$ 4 h^2 (v^1 \J(v^2,v^3) + v^2 \J(v^3,v^1) + v^3 \J(v^1,v^2))
+{\cal O}(h^4)$$
where $\J$ is the Jacobian.
Taking $v^3=1$, for example, recovers the Arakawa Jacobian
Fig. \ref{fig2}(d), and shows that the Arakawa Jacobian also has the
(Casimir) integral $\sum_i u_i$.

It also suggests that in three dimensions with three integrals,
a face-centered-cubic grid (the red points
in a red-black colouring of a cubic grid) is suitable. 
Each vertex is surrounded by 8 tetrahedra. Taking
$\l=((0,0,0),(0,1,1),(1,0,1),(1,1,0))$ (i.e, coupling the
unknowns around one of the tetrahedra) leads to
a fully symmetric discretization of the three-dimensional
Jacobian $\det(\partial v^i/\partial x_j)$.
Using a cubic grid with $\l=((0,0,0),(0,0,1),(0,1,0),(1,0,0))$
leads to a 3D Jacobian with twice the complexity.

\section{Discussion}

We have presented a systematic method for discretizing
PDEs with a known list of integrals. Since all vector
fields $f_i$ with integrals $I^1,\dots,I^p$ can be
written in the form (\ref{Kijk}), the vector fields
$F(\i)$ span all integral-preserving discretizations.
The required symmetry properties of $K$ make the
the finite differences unavoidably complicated,
but sometimes the (optional) spatial symmetries $G$
coincide with the (compulsory) skew symmetries $S_{p+1}$,
reducing the overall complexity of the finite difference.

We close with some comments on future directions.

\begin{enumerate}
\item We have not yet mentioned time integration.
It may not as crucial to preserve integrals in
time as in space; this is not usually done with
the Arakawa Jacobian, for example. If it is
important, we note that linear integrals are
preserved by an consistent linearly covariant method
(such as the Euler method used in (\ref{laxwendroff}));
quadratic integrals are preserved by some Runge-Kutta
methods such as the midpoint rule; and any number
of arbitrary integrals can be preserved by a discrete-time
analogue of (\ref{Kbracket}) \cite{mc-qu-ro}.
With one integral, a simple method is based on
splitting $K$ \cite{mc-qu}.
\item To get simpler finite differences, some of the
spatial symmetries can be broken, as for example in
the half- and quarter-size Arakawa Jacobians in 
Fig. \ref{fig3}(b,c). How important is this in
practice? These Jacobians are still fully translation
invariant. Breaking {\it this} symmetry gives
even simpler tensors $K$, in which, e.g.,
different differences are applied to red and to
black points. Is this useful?
\item Such broken symmetries may be partially
repaired ``on the fly'' during the time integration.
At the $n$th time step we use the finite difference
tensor $\sgn(g_n)gK$, with $g_n$ ranging over the
symmetries. This decreases the symmetry errors
by one power of the time step \cite{is-mc-za}, 
which, with $\Delta t = (\Delta x)^r$, may be plenty.
Most drastically, $K$ could be only first-order
accurate, improving to second through the time
integration. This would require a careful stability
analysis.
\item Although our discretizations are not 
Hamiltonian (unless $\K$ is constant), they can
be volume-preserving. The system (\ref{Kijk}) is
volume preserving for all $I^j$ if $\sum_i \partial K_{ijk\dots}
/\partial u_i=0$ for all $j$, $k,\dots$. In simple
cases this is simple to arrange; $K$ in Eq. (\ref{nonskew})
is not volume-preserving, but $K$ in Eq. 
(\ref{vp}) is. Incorporating volume-preservation
in general is more difficult; see \cite{mc-qu} for
a discussion.
\item We have deliberately avoid mentioning boundaries
and the precise degree of smoothness required of
the arguments that make $\K$ skew. These are studied
in \cite{mc-ro}. If the PDE develops shocks, a
careful weighting of the $F(\i)$ will be required to
capture them well, the analogue of the many methods
for choosing $H$ in (\ref{laxwendroff}) \cite{iserles}.
The present work
applies to infinite or (trivially) to periodic domains.
With finite domains, one can start with $K_\si$ at
an interior point $\i$, and extend it to the boundary
by skew-symmetry, giving finite difference tensors
satisfying certain ``natural'' boundary conditions.
Also on finite domains, there is the possibility of
using global (e.g. spectral) methods. Of course,
these are in the span of our basis, but that is
not the best way to view them. Preserving
integrals with global methods is studied in \cite{mc-ro}.
\item We have concentrated on constructing skew-symmetric
tensors approximating skew operators. Exactly the
same technique can be used to construct symmetric
tensors. We replace the canonical sign function
on $S_{p+1}$ by any sign function $\sigma$ that makes
$S_{p+1}\cong \Z_2$. If $\sigma(\i)=1$ for all $\i$,
for example, the resulting $K$ is completely symmetric,
and when contracted against any $p-1$ of the $I^j$, has
real eigenvalues. If negative definite, the $I^j$ decrease
in time. What is the relationship with the
support operator method \cite{shashkov}?
\end{enumerate}

\subsection*{Acknowledgements}
This paper had a long gestation, during which the support
of the Isaac Newton Institute, Cambridge, and the MSRI,
Berkeley, were invaluable. Useful discussions with
Phil Morrison, Reinout Quispel, Nicolas Robidoux, and
Rick Salmon are gratefully acknowledged.
This work was supported in part by a grant from the
Marsden Fund of the Royal Society of New Zealand.

\section*{Appendix}

Consider the operator $\K=\partial_x$. We were
puzzled by the following: the tensor $K$ in Eq. (\ref{central}),
$$ K = {1\over 2h}\left(\begin{array}{cccc}
\ddots & & & \\
-1 & 0 & 1 & \\
& -1 & 0 & 1 \\
& & & \ddots \\
\end{array}
\right),$$
preserves not just the integral $H$ it operates on,
but also $C=\sum_i u_i$, because $C$ is a Casimir of $K$.
But this two-integral discretization does not arise from any of
the rank 3 skew-tensors we derived in Section 5, Case 2---%
Eq. (\ref{p2d1}) in particular. Contracting with the required
integral $C$ gives a discretization of $\partial_{xxx}$, not
of $\partial_x$. The same is true for any other basis element.

To force $C$ to appear explicitly in the discretization,
we first find a skew differential operator $\J(u,v,w)$ such that
$\J(u,v,\delta C)=\K(u,v)$ for all $\H$. 
If we restrict to a finite domain $D$ so that $\int_D 1\, dx$ is finite,
a natural solution is
$$\J(u,v,w) = u v_x \int_D w\,dx - u w_x \int_D v\, dx + u \int_D 
v w_x\, dx.$$
This is a non-local differential operator, which is the
resolution of the paradox. It only reduces to a local
operator when $w\equiv1$. Discretizing its derivatives by
central derivatives, and integrals $\int w\,dx$  by $\sum_i w_i$,
gives a non-local skew 3-tensor $J$ such that $J(\nabla C)=K$.
In Section 5 we only looked at local tensors.

It seems unlikely that the telescoping sum which makes
this example work will work for nonlinear Casimirs. On the
other hand, it is quite hard to destroy linear ones. Therefore
we suggest the following strategy: temporarily disregard any
known linear integrals (mass, momentum etc.). Construct a skew
tensor so as to preserve the desired nonlinear integrals. Then,
check that this tensor has (some discretization of) 
the required linear integrals as Casimirs. 

The situation is analogous to preserving volume, a linear {\it
differential} invariant discussed in Section 5, note 4.

\end{document}